# BINARY SOLUTIONS FOR OVERDETERMINED SYSTEMS OF LINEAR EQUATIONS

SUBHENDU DAS, CCSI, CALIFORNIA, USA

**Abstract:**

This paper presents a finite step method for computing the binary solution to an overdetermined system of linear algebraic equations $Ax = b$, where $A$ is an m x n real matrix of rank $n < m$, and $b$ is a real m-vector. The method uses the optimal policy of dynamic programming along with the branch and bound concept. Numerical examples are given.

**Introduction**

Many problems in science, engineering, business, management, and economics are formulated as a system of linear equations. In addition many of them also look for the binary or the zero-one solution of these equations. In this paper we solve the following Binary Programming (BP) problem:

$$\text{BP:} \quad Ax = b, \ A \in R^{m \times n}, \ m > n, \ x \in \{0,1\}^n, \ b \in span(A), \ rank(A) = n \tag{1}$$

We define span(A) as the space spanned by the zero-one combinations of columns of A. More precisely

$$span(A) = \{y: \ y = Ax, \ \forall x \in \{0,1\}^n\} \tag{2}$$

It is assumed that the elements of the matrix A are precisely known and the elements of the vector b may have some noise errors. This situation happens in digital communication systems where b is the vector received from the transmitter but the matrix A will be available at both the transmitter and the receiver stations. A detailed description of such a problem can be found in Das [4]. An interesting genesis is provided in Donoho [5].

First we introduce some notations, and then describe the algorithm. Next we show how we create the test problems so that the matrices have the correct ranks and (1) has known solutions. We give two examples. The first one is very small, just enough to show the algorithm details. The second one is little larger and has a noisy b-vector. Then we provide the solution tables for our algorithm. We conclude with discussions of the literature related to our problem.

**Notations**

The column vector *[x₁, x₂, ... , xₙ]'* and its components *{ x₁, x₂, ... , xₙ }* will be denoted by the lower case symbol x. The columns of the matrix A will be similarly represented by *{a₁, a₂, ... , aₙ }*. The i-th state of the system is defined as

$$s_i = \begin{cases} s_{i-1}, & \text{if } x_i = 0 \\ s_{i-1} - a_i & \text{if } x_i = 1 \end{cases}, \ s_0 = b, \ i = 1, \dots, n-1 \tag{3}$$

As an example, when the decision *x₁ = 0* is used then the state is *s₁ = b* and when the decision *x₁ = 1* is chosen then the state is *s₁ = b - a₁*. The state is related to the right hand side of (1). When a column is removed from the matrix A, it goes to the right hand side of (1) multiplied by the corresponding value of the variable.

**Algorithm**

One of the major concepts we use is the Bellman's dynamic programming (DP) principle of optimality. This concept can be found in many text books on operations research, Wolsey [12] or dynamic programming, Nemhauser [10]. The DP principle is stated in the following way – "An optimal set of decisions has the property that whatever the first decision is, the remaining decisions must be optimal with respect to the outcome which results from the first decision".

Thus when we look for optimal value for $x_1$ we first choose $x_1 = 0$ and then find the optimal values of all other variables that give the best solution for $B_1 x = b$ with $x$ as unconstrained to any real number. Here $B_1$ is the $A$ matrix, with the first column removed, and $x$ represents the remaining variables. This unconstrained problem is solved using the pseudo inverse, see Golub[4] (p. 257). This is the minimum summed squared error (SSE) or the least square solution for the state $s_1 = b$, which is the right hand side of the equation (1). Next we choose $x_1 = 1$, and use the same method to get the SSE for the state $s_1 = b - a_1$ of the system. The optimal decision for $x_1$ is the minimum of the two SSE results. The variable $x_1$ and the column $a_1$ are moved to the right hand side of (1). The new state, $s_1$, becomes the optimal state for the next variable, $x_2$. Thus the foundation of the algorithm is clearly based on the DP principle of optimality.

The branch and bound (BB) method is embedded in the above algorithm, because we are using the relaxed problem, that is the unconstrained problem, to select the bounds for the two SSE values. Note that when we are deciding for $x_3$ variable for example, there is no need to consider all possibilities for $x_1$ and $x_2$ again. That is because we have considered all possible values for $x_3$ when we considered previous variables. When we were deciding for $x_1$ we did indeed consider all possible combinations of $x_2$ and $x_3$ using the relaxation logic of pseudo inverse.

The formal steps of the procedure can be written as:

**DPBB Algorithm:** (4)

For i = 1 to n - 1 repeat

    Select $x_i = 0$, state $s_i = s_{i-1}$, and $B_i = [\,a_{i+1},\ldots,a_n]$

    Find pseudo inverse $P_i$ of $B_i$

    $x_{i0} = P_i * s_i$                           //optimal values for the remaining variables

    $s_{i0} = B_i * x_{i0}$                         //optimal estimate of the state using remaining variables

    $SSE_{i0} = (s_i - s_{i0})'\,(s_i - s_{i0})$     //state estimation error

    Select $x_i = 1$, state $s_i = s_{i-1} - a_i$,

    $x_{i1} = P_i * s_i$                            //same as above

    $s_{i1} = B_i * x_{i1}$

    $SSE_{i1} = (s_i - s_{i1})'\,(s_i - s_{i1})$

    If $SSE_{i1} < SSE_{i0}$ select $x_i = 1$ Else $x_i = 0$

End of For loop

In the above algorithm $x_{i0}$ and $x_{i1}$ are the unconstrained optimal values of the remaining x variables. These x variables are then used to generate the best estimate of the current state $s_{i0}$ and $s_{i1}$. The decision for the final variable $x_n$ is similar. We use a different method because there are no remaining columns of A. This part of the algorithm can be written in the following way:

**DPBB Algorithm – Last:** (5)

For i=n do the following:

    Select state $s_n = s_{n-1}$

    Select $x_n = 0$
    $SSE_{n0} = s_n's_n$         //the error is the magnitude of the state

    Select $x_n = 1$
    $SSE_{n1} = (s_n-a_n)'(s_n-a_n)$     //the error is the magnitude of the state

    If $SSE_{n1} < SSE_{n0}$ select $x_n = 1$ Else $x_n = 0$

End of the trivial For loop.

It should be clear that our problem is essentially the same as the least square solution problem of the standard type defined by:

$$\min_{x \in \{0,1\}^n}\{\|Ax - b\|_2\}$$

In our computer program we just did not take the square root of the 2-norm and used SSE instead. Observe that the over determined system usually does not have exact solution, Golub [7, p 236]. Since the matrix A is of full rank the value for SSE can never be zero also.

**Problem Construction**

In this section we show how we have constructed the two test examples for demonstrating our algorithm. The first problem has three unknowns and has ten equations. The second problem has ten unknowns and has twenty equations. After we describe the problem constructions method, we will walk through the first example to illustrate the algorithm in details and then simply give the partial results for the second example. Both problems are constructed using the same principle. We have used Mathematica software tools for our analysis.

We represent the data in the columns of the A matrix as the digital sample values of several independent functions. This process ensures that the columns are independent and therefore has full column rank. For the first problem the following three functions were used, since there are three unknowns:

$T = 0.001, \quad f_0 = \frac{1}{T}$

$g_1(t) = \text{Sin}(2\pi f_0 t) + 1,$  (6)

$g_2(t) = \text{Cos}(2\pi f_0 t) + 1, \quad g_3(t) = \text{Cos}(4\pi f_0 t) + 1$

Ten samples were generated using the equal sample interval of *dT = T/10*. The first sample started at *t = dT*.

For the second problem we use the following 10 functions, one for each column of the matrix A. The remaining parameters were same as problem one. This problem has both negative and positive elements in the column vectors.

$h_1(t) = \text{Sin}(2\pi f_0 t) \qquad h_2(t) = \text{Cos}(2\pi f_0 t) \quad h_3(t) = \text{Cos}(4\pi f_0 t)$

$h_4(t) = \text{Sin}(4\pi f_0 t) \qquad h_5(t) = \text{Sin}(\pi f_0 t) \quad h_6(t) = \text{Cos}(\pi f_0 t)$  (7)

$h_7(t) = \text{Cos}(6\pi f_0 t) \qquad h_8(t) = \text{Sin}(6\pi f_0 t)$

$h_9(t) = 1 - e^{-3000t} \qquad h_{10}(t) = e^{-3000t}$

For generating the b vector we selected x = {1, 0, 1} for problem 1. Therefore our algorithm should produce the above x values as the correct solution. For the second problem we selected arbitrarily x = {1,0,1,1,1,0,1,1,1,0}, and therefore this is the correct solution too.

For the second problem we also used a uniform random number generator for generating twenty random numbers between 0 and 1. Then 20% of these numbers, *w = {$w_i$, i=1,...,20}*, were added to each element of the b-vector to introduce some random noise. It is assumed that this will be the residual noise in the b-vector after it has been processed using other signal processing algorithms, like the finite impulse response filters (FIR), see Lyons [8], on the original noisy vector. For the first problem no noise was considered.

**Example 1**

As stated before this example has three unknowns and ten equations. Using the matrix notation we present the given data as follows:

$$\begin{bmatrix} 1.4683 \\ 1.0955 \\ 1.0855 \\ 1.4683 \\ 1.5125 \\ 0.86579 \\ 0.12058 \\ 0.12058 \\ 0.86579 \\ 1.5125 \end{bmatrix} = \begin{bmatrix} 0.81381 & 0.90451 & 0.65451 \\ 1 & 0.65451 & 0.095492 \\ 1 & 0.34549 & 0.085492 \\ 0.81381 & 0.095492 & 0.65451 \\ 0.51254 & 0 & 1 \\ 0.21128 & 0.095492 & 0.65451 \\ 0.025086 & 0.34549 & 0.095492 \\ 0.025086 & 0.65451 & 0.095492 \\ 0.21128 & 0.90451 & 0.65451 \\ 0.51254 & 1 & 1 \end{bmatrix} \begin{bmatrix} x_1 \\ x_2 \\ x_3 \end{bmatrix}$$

Our problem is to find out the *0-1* solution of the above linear system of equations for the unknown components of the vector $x = \{ x_1, x_2, x_3 \}$.

The DPBB algorithm (4) selects the unknown variables one by one, in sequence, and starting from the first variable $x_1$. When we select the first variable $x_1$ for decisions, we take the last two columns of the matrix $A$ and call it matrix $B_1$. The pseudo inverse $P_1$ of the matrix $B_1$ is given by the formula, see Golub [7] on p. 257:

$$P_1 = (B_1' B_1)^{-1} B_1' \qquad (8)$$

For this example the matrices $B_1$ and $P_1$ are given below

$$B_1 = \begin{bmatrix} 0.90451 & 0.65451 \\ 0.65451 & 0.095492 \\ 0.34549 & 0.085492 \\ 0.095492 & 0.65451 \\ 0 & 1 \\ 0.095492 & 0.65451 \\ 0.34549 & 0.095492 \\ 0.65451 & 0.095492 \\ 0.90451 & 0.65451 \\ 1 & 1 \end{bmatrix} \qquad P_1 = \begin{bmatrix} 0.22472 & 0.28360 & 0.13527 & -0.16360 & -0.32 & -0.16360 & 0.13527 & 0.28360 & 0.22472 & 0.16 \\ 0.02472 & -0.16360 & -0.06472 & 0.28360 & 0.47999 & 0.28360 & -0.06472 & -0.16360 & 0.02472 & 0.16 \end{bmatrix}$$

For the selection $x_1 = 0$, the state is $s_1$, and is equal to the *b*-vector. We multiply $P_1$ by $s_1$ to get the optimal values of the remaining *x* variables $\{x_2, x_3\}$. Then we multiply $B_1$ by these x variables to get the optimal estimate for the state $s_1$. The difference between $s_1$ and its optimal estimate will produce the SSE for this selection of $x_1$. For $x_1 = 0$ selection, we have $[x_2 \; x_3]' = P * s_1 = P*b = [0.410034, 1.41003]'$. Then we use this unconstrained optimal values for $x = [x_2 \; x_3]'$ to get the estimate for $s_1$ using $B_1 * x$. The tables below show the state, the optimal estimate of the state, the error values, and the SSE. The tables are identified as 1A and 1B for the two selections, $x_1 = 0$ and $x_1 = 1$, respectively.

| Decision variable is $x_1$. | Decision is $x_1 = 0$. | | $B = [a_2, a_3]$ is used for pseudo inverse. | | | Initial state $s_0 = b$ | | | |
|---|---|---|---|---|---|---|---|---|---|
| Optimal choice for remaining variables: X2= 0.410034     X3=1.41003 | | | | | | | | | |
| Initial state $s_0$ | 1.4683 | 1.0954 | 1.0854 | 1.4683 | 1.5125 | 0.8657 | 0.1205 | 0.1250 | 0.8657 | 1.5125 |
| New state $s_1 = s_0$ | 1.4683 | 1.0954 | 1.0854 | 1.4683 | 1.5125 | 0.8657 | 0.1205 | 0.1250 | 0.8657 | 1.5125 |
| State estimate | 1.2937 | 0.4030 | 0.2763 | 0.9620 | 1.4100 | 0.9620 | 0.2763 | 0.4030 | 1.2937 | 1.8200 |
| Estimation error | 0.1745 | 0.6924 | 0.8191 | 0.5062 | 0.1025 | -0.0962 | -0.1557 | -0.2824 | -0.4279 | -0.3075 |
| SSE | 1.8389 | Table 1A – Problem 1 | | | | | | | | |

Now repeat the above procedure for the selection $x_1 = 1$ for the same variable $x_1$. The state vector is $s_1 = b - a_1$.

| Decision variable is $x_1$. | Decision is $x_1 = 1$. | | $B = [a_2, a_3]$ is used for pseudo inverse. | | | Initial state $s_0 = b$ | | | |
|---|---|---|---|---|---|---|---|---|---|
| Optimal choice for remaining variables: X2 = 0    X3 = 1.0 | | | | | | | | | |
| Initial state $s_0$ | 1.4683 | 1.0954 | 1.0854 | 1.4683 | 1.5125 | 0.8657 | 0.1205 | 0.1250 | 0.8657 | 1.5125 |
| $a_1$ | 0.8138 | 1.0 | 1.0 | 0.8138 | 0.5125 | 0.2112 | 0.0250 | 0.0250 | 0.2112 | 0.5125 |
| New state $s_1 = s_0 - a_1$ | 0.6545 | 0.0954 | 0.0855 | 0.6545 | 1.0 | 0.6545 | 0.0954 | 0.0954 | 0.6545 | 1.0 |
| State estimate | 0.6545 | 0.9545 | 0.9845 | 0.6545 | 1.0 | 0.6545 | 0.0954 | 0.0954 | 0.6545 | 1.0 |
| Estimation error | 1.1E-16 | -6.9E-17 | -6.9E-17 | 1.1E-16 | 1.1E-16 | 1.1E-16 | 1.3E-17 | 1.3E-17 | 1.1E-16 | 1.1E-16 |
| SSE | 8.39E-32 | Table 1B – Problem 1 | | | | | | | | |

Since the SSE is lower for $x_1 = 1$, we decide that the optimal value for $x_1$ is $1$ and that is the correct result as used in the problem formulation step.

The decision tables, 2A and 2B, for the second variable $x_2$ are similarly computed and shown below. In this case the starting state is $s_2 = b - a_1$ since $x_1$ was found as $1$. The matrix $B_2$ now is the last column of the matrix A. The pseudo inverse $P_2$ in this case is given by

$P_2$ = [0.174536, 0.0254645, 0.0254645, 0.174536, 0.266666, 0.174536, 0.0254645, 0.0254645, 0.174536, 0.266666]

The corresponding unconstrained optimal value for $x_3$ is 1.0. Note that we are now working for $x_2$ variable. For $x_2 = 0$ we get the following table

| Decision variable is $x_2$. | Decision is $x_2 = 0$. | | $B = [a_3]$ is used for pseudo inverse. | | | Initial state $s_1 = b - a_1$ | | | |
|---|---|---|---|---|---|---|---|---|---|
| Optimal choice for remaining variables: X3=1.0 | | | | | | | | | |
| Initial state $s_1$ | 0.6545 | 0.0954 | 0.0854 | 0.6545 | 1.0 | 0.6545 | 0.0954 | 0.0954 | 0.6545 | 1.0 |
| New state $s_2 = s_1$ | 0.6545 | 0.0954 | 0.0854 | 0.6545 | 1.0 | 0.6545 | 0.0954 | 0.0954 | 0.6545 | 1.0 |
| State estimate | 0.6545 | 0.0954 | 0.0854 | 0.6545 | 1.0 | 0.6545 | 0.0954 | 0.0954 | 0.6545 | 1.0 |
| Estimation error | 1.11E-16 | -6.93E-17 | -6.93E-17 | 1.11E-16 | 1.11E-16 | 1.11E-16 | 1.38E-17 | 1.38E-17 | 1.11E-16 | 1.11E-16 |
| SSE | 8.39E-32 | Table 2A – Problem 1 | | | | | | | | |

For the selection $x_2 = 1$ we create a similar table, along the lines of the dynamic programming theory.

| Decision variable is $x_2$. | Decision is $x_2 = 1$. | | B = [$a_3$] is used for pseudo inverse. | | | Initial state $s_1 = b - a_1$ | | | | |
|---|---|---|---|---|---|---|---|---|---|---|
| Optimal choice for remaining variables: X3=0.33333 | | | | | | | | | | |
| Initial state $s_1$ | 0.6545 | 0.0954 | 0.0854 | 0.6545 | 1.0 | 0.6545 | 0.0954 | 0.0954 | 0.6545 | 1.0 |
| $a_2$ | 0.9045 | 0.6545 | 0.3454 | 0.0954 | 0.0 | 0.0954 | 0.3454 | 0.6545 | 0.9045 | 1.0 |
| New state $s_2 = s_1 - a_2$ | -0.25 | -0.5590 | -0.2499 | 0.5590 | 1.0 | 0.5590 | -0.2499 | -0.5590 | -0.25 | 0.0 |
| State estimate | 0.2181 | 0.0318 | 0.0318 | 0.2181 | 0.3333 | 0.2181 | 0.0318 | 0.0318 | 0.2181 | 0.3333 |
| Estimation error | -0.4681 | -0.5908 | -0.2818 | 0.3408 | 0.6666 | 0.3408 | -0.2818 | -0.5908 | -0.4681 | -0.3333 |
| SSE | 2.0833 | Table 2B – Problem 1 | | | | | | | | |

The optimal decision for $x_2$ is $x_2 = 0$, since the SSE for the first table, 2A, is lower and that is the correct choice also as defined in the formulation stage of the problem.

In the last step, for the variable $x_3$, we do not use the pseudo inverse, see (5) for the DPBB algorithm. The state is still $b-a_1$ because the optimal value for $x_2$ turned out to be zero. For the selection $x_3 = 0$, SSE is computed using the estimation error, which in this case is just the magnitude of the new state which is $b-a_1$.

| Decision variable is $x_3$. | Decision is $x_3 = 0$. | | B is not used for pseudo inverse. | | | Initial state $s_2 = s_1 = b - a_1$ | | | | |
|---|---|---|---|---|---|---|---|---|---|---|
| Initial state $s_2$ | 0.6545 | 0.0954 | 0.0854 | 0.6545 | 1.0 | 0.6545 | 0.0954 | 0.0954 | 0.6545 | 1.0 |
| New state $s_3 = s_2$ | 0.6545 | 0.0954 | 0.0854 | 0.6545 | 1.0 | 0.6545 | 0.0954 | 0.0954 | 0.6545 | 1.0 |
| Estimation error | 0.6545 | 0.0954 | 0.0854 | 0.6545 | 1.0 | 0.6545 | 0.0954 | 0.0954 | 0.6545 | 1.0 |
| SSE | 3.7500 | Table 3A – Problem 1 | | | | | | | | |

For the choice of $x_3 = 1$, we generate the following table, 3B, using the same method. The new state is $s_2 - a_2$.

| Decision variable is $x_3$. | Decision is $x_3 = 1$. | | B is not used for pseudo inverse. | | | Initial state $s_2 = s_1 = b - a_1$ | | | | |
|---|---|---|---|---|---|---|---|---|---|---|
| Initial state $s_2$ | 0.6545 | 0.0954 | 0.0854 | 0.6545 | 1.0 | 0.6545 | 0.0954 | 0.0954 | 0.6545 | 1.0 |
| $a_3$ | 0.6545 | 0.0954 | 0.0854 | 0.6545 | 1.0 | 0.6545 | 0.0954 | 0.0954 | 0.6545 | 1.0 |
| New state $s_3 = s_2 - a_3$ | 0. | -8.32E-17 | -8.32E-17 | 0 | 0 | 0 | 0 | 0 | 0 | 0 |
| Estimation error | 0. | -8.32E-17 | -8.32E-17 | 0 | 0 | 0 | 0 | 0 | 0 | 0 |
| SSE | 1.38E-32 | Table 3B – Problem 1 | | | | | | | | |

Thus the optimal decision for $x_3$ is *1* since SSE for the second table, 3B, is lower and is also the correct decision. This concludes the implementation and verification of the algorithm for the example problem one. We have found the correct optimal solution in three steps for the three unknown variables of the problem.

**Example 2**

We briefly describe the solution for the sample problem two. This is also a small size problem but with a noisy b-vector. This problem has 10 unknown variables and 20 equations. The columns of the matrix A are generated from the samples of the h functions defined in (7). We do not do any filtering in this example, instead only model the filtered vector b by adding some residual noise *w*. The tables are provided for the numerically oriented readers.

The data for the matrix A is given below. To save space the complete data is not provided. It only shows some elements of the data and hopes to provide enough information so that any one will be able to reproduce the same

| A matrix columns |||||||||| b Vectors ||
| 1 | 2 | 3 | 4 | 5 | 6 | 7 | 8 | 9 | 10 | b | b + w |
|---|---|---|---|---|---|---|---|---|---|---|---|
| 0.2955 | 0.9663 | 0.8355 | 0.5649 | 0.1494 | 1. | 0.6234 | 0.8019 | 0.1412 | 1. | 3.4122 | 3.5466 |
| 0.5649 | 0.8355 | 0.3694 | 0.9334 | 0.2955 | 0.9663 | -0.2225 | 1. | 0.2636 | 0.8668 | 3.2046 | 3.2808 |
| 0.7840 | 0.6305 | -0.2250 | 0.9776 | 0.4351 | 0.9111 | -0.9009 | 0.4450 | 0.3698 | 0.7514 | 1.8856 | 2.0184 |
| -------- | -------- | -------- | -------- | -------- | -------- | -------- | -------- | -------- | -------- | -------- | -------- |
| -0.7840 | 0.6305 | -0.2250 | -0.9776 | 0.4351 | -0.9111 | -0.9009 | -0.4450 | 0.9798 | 0.0881 | -1.9178 | -1.8409 |
| -0.5649 | 0.8355 | 0.3694 | -0.9334 | 0.2955 | -0.9663 | -0.2225 | -1. | 0.9906 | 0.0764 | -1.0652 | -0.9618 |
| -0.2955 | 0.9663 | 0.8355 | -0.5649 | 0.1494 | -1. | 0.6234 | -0.8019 | 1. | 0.0662 | 0.9461 | 1.1094 |

data and the same results if desired. The original b vector corresponding to the correct *x* vector, and the noisy vector, *b + w*, are also in the same table.

We also do not give all the tables for solving this entire problem. It will require 10 tables, one for each variable; each table will have 20 columns, which will be too big for the space allowed for this paper, and probably is not necessary also. Thus we give the tables for the two decision options of the first variable $x_1$, just to identify the nature of the tables and the associated computer data structure for a larger problem. We also do not give all the data in each table; some columns are removed to fit the table in the page with a readable font size. The first table for the selection $x_1 = 0$ is shown in 1A.

| Decision variable is $x_1$. | Decision is $x_1 = 0$. ||| B = [$a_2, a_3,…,a_{10}$] is used for pseudo inverse. |||| Initial state $s_0$ = b + w |||
|---|---|---|---|---|---|---|---|---|---|---|
| Optimal choice for remaining variables: -1.1691, 0.8434, 0.9172,-3.8932, 3.3059, 0.9510, 1.0410, 6.3035, -1.6104 ||||||||||||
| Initial state $s_0$ | 3.5466 | 3.2808 | 2.0184 | 0.7165 | 0.2800 | | 0.4550 | -1.0923 | -1.8408 | -0.9618 | 1.1094 |
| New state $s_1 = s_0$ | 3.5466 | 3.2808 | 2.0184 | 0.7165 | 0.2800 | | 0.4550 | -1.0923 | -1.8408 | -0.9618 | 1.1094 |
| State estimate | 3.5247 | 3.3303 | 2.0153 | 0.6732 | 0.2523 | | 0.4727 | -1.1037 | -1.8154 | -0.9981 | 1.1239 |
| Estimation error | 0.0218 | -0.0495 | 0.0030 | 0.0433 | 0.0277 | | -0.0177 | 0.0113 | -0.0253 | 0.0363 | -0.0145 |
| SSE | 0.0318 | Table 1A – Problem 2 ||||||||||

The second table for the decision $x_1 = 1$ is created in a similar way and is shown in table 1B.

| Decision variable is $x_1$. | Decision is $x_1 = 1$. | | | B = [$a_2, a_3,...,a_{10}$] is used for pseudo inverse. | | Initial state $s_0$ = b + w | | | |
|---|---|---|---|---|---|---|---|---|---|
| Optimal choice for remaining variables:  -0.7086, 0.8640, 1.0139, -0.9105, 0.1759, 0.9533, 1.0600, 2.4674, 0.7896 | | | | | | | | | |
| Initial state $s_0$ | 3.5466 | 3.2808 | 2.0184 | 0.7165 | 0.280 | -1.0923 | -1.8408 | -0.9618 | 1.1094 |
| $a_1$ | 0.2955 | 0.5649 | 0.7840 | 0.9334 | 1.0 | -0.9334 | -0.7840 | -0.5649 | -0.2955 |
| New state $s_1 = s_0 - a_1$ | 3.2510 | 2.7159 | 1.2344 | -0.2168 | -0.7199 | -0.1588 | -1.0568 | -0.3969 | 1.4050 |
| State estimate | 3.2325 | 2.7576 | 1.2329 | -0.2549 | -0.7466 | -0.1752 | -1.0325 | -0.4261 | 1.4163 |
| Estimation error | 0.0184 | -0.0417 | 0.0014 | 0.0380 | 0.0267 | 0.0163 | -0.0243 | 0.0291 | -0.0113 |
| SSE | 0.0291 | Table 1B – Problem 2 | | | | | | | |

From the above two tables we see that SSE is lower in the table for $x_1 = 1$ and therefore the optimal decision for $x_1$ is *1* which is also the correct result as defined during the problem formulation stage. We now give, just for the sake of a feeling of completeness, the last two tables for the last decision variable $x_{10}$.

| Decision variable is $x_{10}$. | Decision is $x_{10} = 0$. | | | B is not used for pseudo inverse. | | Initial state $s_9$ | | | |
|---|---|---|---|---|---|---|---|---|---|
| Initial state $s_9$ | 0.1344 | 0.0762 | 0.1328 | 0.1616 | 0.1294 | 0.1919 | 0.1881 | 0.0769 | 0.1033 | 0.1633 |
| New state $s_{10} = s_9$ | 0.1344 | 0.0762 | 0.1328 | 0.1616 | 0.1294 | 0.1919 | 0.1881 | 0.0769 | 0.1033 | 0.1633 |
| Estimation error | 0.1344 | 0.0762 | 0.1328 | 0.1616 | 0.1294 | 0.1919 | 0.1881 | 0.0769 | 0.1033 | 0.1633 |
| SSE | 0.3763 | Table 2A –Problem 2 | | | | | | | |

Note that for the last stage there is no pseudo inverse. The estimate is the magnitude of the final state. So the error is also based on the final state by default.

| Decision variable is $x_{10}$. | Decision is $x_{10} = 1$. | | | B is not used for pseudo inverse. | | Initial state $s_9$ | | | |
|---|---|---|---|---|---|---|---|---|---|
| Initial state $s_9$ | 0.1344 | 0.0762 | 0.1328 | 0.1616 | 0.1294 | 0.1919 | 0.1881 | 0.0769 | 0.1033 | 0.1633 |
| $a_{10}$ | 1 | 0.8668 | 0.7514 | 0.6514 | 0.5647 | 0.1173 | 0.1017 | 0.0881 | 0.0764 | 0.0662 |
| New state $s_{10} = s_9 - a_{10}$ | -0.8655 | -0.7906 | -0.6186 | -0.4897 | -0.4352 | 0.0746 | 0.0864 | -0.0112 | 0.0269 | 0.0970 |
| Estimation error | -0.8655 | -0.7906 | -0.6186 | -0.4897 | -0.4352 | 0.0746 | 0.0864 | -0.0112 | 0.0269 | 0.0970 |
| SSE | 2.67086 | Table 2B – Problem 2 | | | | | | | |

Since the lower SSE is obtained from the first table, 2A, the optimal solution for state $x_{10}$ is *0*, which also is the correct solution.

We have used pseudo inverse (7) of rectangular matrices in all our DPBB algorithm steps. It is interesting to point out that a direct application of pseudo inverse to the entire problem will also give correct result in the absence of any noise in the b vector of the system. As an example the second problem gives the following result for x as shown in the table below with the direct application of the formula (9) to the matrix A of (1):

$$x = P.b, \quad where \ \ P = [A'A]^{-1}A \tag{9}$$

| Direct application of Pseudo Inverse for the b vector | | | | | | | | | |
|---|---|---|---|---|---|---|---|---|---|
| x | 1.0 | 0.E-3 | 1.000 | 1.00 | 1. | 0.E-2 | 1.0000 | 1.0000 | 1.0 | 0.E-2 |
| True solution | 1 | 0 | 1 | 1 | 1 | 0 | 1 | 1 | 1 | 0 |

However if we use the noisy $b + w$ vector in the above formula (9) then the pseudo inverse gives the following normalized solution:

| Direct application of Pseudo Inverse for the noisy b + w vector | | | | | | | | | |
|---|---|---|---|---|---|---|---|---|---|
| x | 0.368 | 0.092 | 0.063 | 0.096 | 0.843 | -0.936 | 0.063 | 0.075 | -1.0 | 0.779 |
| True solution | 1 | 0 | 1 | 1 | 1 | 0 | 1 | 1 | 1 | 0 |

It is clear from the last table that we cannot extract the correct solution from the x values of the previous table corresponding to the noisy b + w vector. Therefore the DPBB algorithm works better in more realistic environment.

**Discussions**

It should be clear that the BP problem defined by (1) and the corresponding DPBB algorithm defined by (4) is not a NP complete problem. However, a problem very similar to (1) has been defined as NP-complete problem by Murty and Kabadi [9]. This NP-complete problem is stated as follows. Given the positive integers $\{d_0, d_1, \ldots, d_n\}$, is there a solution to:

$$\sum_{i=1}^{n} d_i y_i = d_0 \ , \ \ y_i = 0 \text{ or } 1, \ \ for \ all \ i \quad\quad\quad (10)$$

Our problem is very similar to (10); except the numbers in problem (1) are all real numbers. And also (10) is an underdetermined system and (1) is an overdetermined system. The problem (10) has also been listed in Garey [6, p. 223].

Linear Integer Programming (LIP) is a well known approach for solving problems similar to the BP problem defined in (1). For small size problems the LIP approach is very effective. The LIP requires an optimization criterion and we do not have any such objective function with our problem (1). However, it is well known that the linear integer programming is an NP complete problem; see Wolsey [12, p. 103]. The zero-one linear integer programming is also listed as NP complete problem in Garey, [6, p. 245].

It is also possible to solve the BP problem (1) by converting it to a binary quadratic programming problem. YoonAnn [13] has used such an approach for an overdetermined system. A class of binary quadratic programming problem, such as with non-negative coefficients for the quadratic terms, is also known as NP problem; see Garey [6, p. 245]. Also see Axehill [1]. The method discussed in the present paper solves the problem in polynomial time.

There are many numerical methods available for solving the overdetermined linear system of equations. But most of them are for real valued solutions. It seems that the there has not been much on work done for the binary solutions of the problem defined in (1). The literature search did not produce any such numerical work. The general approach seems to convert the problem to an integer optimization problem which in most cases requires NP algorithms.

The numerical error in using the pseudo inverse may be of concern as mentioned in Cline [3]. He suggests some alternative decomposition methods. However we are using Mathematica software tool, which allows calculations with an accuracy of any preselected number of digits. As an example, using this tool, all calculations can be performed with 100 decimal digits of accuracy without any noticeable difference in computational time on a standard laptop computer. Thus numerical error does not appear to be of particular concern for the problem in (1).

If we replace the pseudo inverse algorithm by some other numerically efficient algorithms like Cholesky factorization then each inverse will require $(m+n/3)n^2$ operations as shown in Golub [7, p. 238]. Thus it is clear that DPBB is a polynomial time algorithm and therefore again (1) is not a NP-Complete problem.

Bellman [2] has first shown how the solution problem of a set of linear simultaneous equations, with positive definite square matrix, can be converted to a multistage decision problem using his dynamic programming (DP) principle. Later Roger [11] has shown how this DP principle can be implemented using analytical expressions for the case of the overdetermined systems. In the present paper we have extended this DP principle numerically, together with the BB concept, to binary solutions for the overdetermined systems.

**Conclusions**

We have given a straight forward computational procedure for finding the binary solutions of the overdetermined systems of linear equations. The procedure takes only n steps; where n is the number of unknown variables in the equations. The algorithm is a polynomial time process.